# MINIMAL SUFFICIENT CAUSATION AND DIRECTED ACYCLIC GRAPHS

By Tyler J. VanderWeele and James M. Robins

*University of Chicago and Harvard University*

Notions of minimal sufficient causation are incorporated within the directed acyclic graph causal framework. Doing so allows for the graphical representation of sufficient causes and minimal sufficient causes on causal directed acyclic graphs while maintaining all of the properties of causal directed acyclic graphs. This in turn provides a clear theoretical link between two major conceptualizations of causality: one counterfactual-based and the other based on a more mechanistic understanding of causation. The theory developed can be used to draw conclusions about the sign of the conditional covariances among variables.

**1. Introduction.** Two broad conceptualizations of causality can be discerned in the literature, both within philosophy and within statistics and epidemiology. The first conceptualization may be characterized as giving an account of the effects of certain causes; the approach addresses the question, "Given a particular cause or intervention, what are its effects?" In the contemporary philosophical literature, this approach is most closely associated with Lewis' work [17, 18] on counterfactuals. In the contemporary statistics literature, this first approach is closely associated with the work of Rubin [30, 31] on potential outcomes, of Robins [25, 26] on the use of counterfactual variables in the context of time-varying treatment and of Pearl [21] on the graphical representation of various counterfactual relations on directed acyclic graphs. This counterfactual approach has been used extensively in statistics both in the development of theory and in application. The second conceptualization of causality may be characterized as giving an account of the causes of particular effects; this approach attempts to address the question, "Given a particular effect, what are the various events which might have been its cause?" In the contemporary philosophical literature, this second approach is most notably associated with Mackie's work [19] on insufficient but necessary components of unnecessary but sufficient conditions (INUS conditions) for an effect. In the epidemiologic literature, this approach is most closely associated with Rothman's work [29] on sufficient-component causes. This work is more closely related to the various mechanisms for a particular effect than is the counterfactual approach. Rothman's work on sufficient-component causes has, however, seen relatively little development, extension or application, though the basic framework is routinely taught in introductory epidemiology courses. Perhaps the only major attempt in the statistics literature to extend and apply Rothman's theory has been the work of Aickin [1] (comments relating Aickin's work to the present work are available from the authors upon request).

In this paper, we incorporate notions of minimal sufficient causes, corresponding to Rothman's sufficient-component causes, within the directed acyclic graph causal framework [21]. Doing so essentially unites the mechanistic and the counterfactual approaches into a single framework. As will be seen in Section 5, we can use the framework developed to draw conclusions about the sign of the conditional covariances among variables. Without the theory developed concerning minimal sufficient causes, such conclusions cannot be drawn from causal directed acyclic graphs. In a related paper [35] we have discussed how these ideas relate to epidemiologic research. The present paper develops the theory upon which this epidemiologic discussion relies.

The theory developed in this paper is motivated by several other considerations. As will be seen below, the incorporation of minimal sufficient cause nodes allows for the identification of certain conditional independencies which hold only within a particular stratum of the conditioning variable (i.e., "asymmetric conditional independencies," [7]) which were not evident without the minimal sufficient causation structures. We note that these asymmetric conditional independencies have been represented elsewhere by Bayesian multinets [7] or by trees [3]. Another motivation for the development of the theory in this paper concerns the notion of interaction. Product terms are frequently included in regression models to assess interactions among variables; these statistical interactions, however, even if present, need not imply the existence of an actual mechanism in which two distinct causes both participate. Interactions which do concern the actual mechanisms are sometimes referred to as instances of "synergism" [29], "biologic interactions" [32] or "conjunctive causes" [20], and the development of minimal sufficient cause theory provides a useful framework to characterize mechanistic interactions. In related work [37] we have derived empirical tests for interactions in this sufficient cause sense.

As yet further motivation, we conclude this Introduction by describing how the methods we develop in this paper clarified and helped resolve an analytic puzzle faced by psychiatric epidemiologists. Consider the following somewhat simplified version of a study reported in Hudson et al. [10]. Three hundred pairs of obese siblings living in an ethnically homogenous upper-middle class suburb of Boston are recruited and cross classified by the presence or absence of two psychiatric disorders: manic-depressive disorder $P$ and binge eating disorder $B$. The question of scientific interest is whether these two disorders have a common genetic cause, because, if so, studies to search for a gene or genes that cause both disorders would be useful. Consider two analyses. The first analysis estimates the covariance $\beta$ between $P_{2i}$ and $B_{1i}$, while the second analysis estimates the conditional covariance $\alpha$ between $P_{2i}$ and $B_{1i}$ among subjects with $P_{1i} = 1$, where $B_{ki}$ is 1 if the $k$th sibling in the $i$th family has disorder $B$ and is zero otherwise, with $P_{ki}$ defined analogously. It was found that the estimates $\beta$ and $\alpha$ were both positive with 95% confidence intervals that excluded zero.

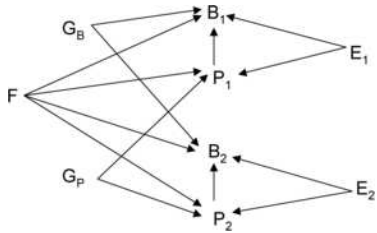

Fig. 1. *Causal directed acyclic graph under the alternative hypothesis of familial coaggregation.*

Hudson et al. [10] substantive prior knowledge is summarized in the directed acyclic graph of Figure 1 in which the $i$ index denoting family is suppressed. In what follows, we will make reference to some standard results concerning directed acyclic graphs; these results are reviewed in detail in the following section.

In Figure 1, $G_B$ and $G_P$ represent the genetic causes of $B$ and $P$, respectively, that are not common causes of both $B$ and $P$. The variables $E_1$ and $E_2$ represent the environmental exposures of siblings 1 and 2, respectively, that are common causes of both diseases, for example, exposure to a particularly stressful school environment. The variables $G_B$ and $G_P$ are assumed independent as would typically be the case if, as is highly likely, they are not genetically linked. Furthermore, as is common in genetic epidemiology, the environmental exposures $E_1$ and $E_2$ are assumed independent of the genetic factors. The causal arrows from $P_1$ to $B_1$ and $P_2$ to $B_2$ represent the investigators' beliefs that manic-depressive disorder may be a cause of binge eating disorder but not vice-versa. The node $F$ represents the common genetic causes of both $P$ and $B$ as well as any environmental causes of both $P$ and $B$ that are correlated within families. There is no data available for $G_B$, $G_P$, $E_1$, $E_2$ or $F$. The reason for grouping the common genetic causes with the correlated environmental causes in $F$ is that, based on the available data $\{P_{ki}, B_{ki};\ i = 1, \ldots, 300, k = 1, 2\}$, we can only hope to test the null hypothesis that $F$ so defined is absent, which is referred to as the hypothesis of no familial coaggregation. If this null hypothesis is rejected, we cannot determine from the available data whether $F$ is present due to a common genetic cause or a correlated common environmental cause. Thus $E_1$ and $E_2$ are independent on the graph because, by definition, they represent the environmental common causes of $B$ and $P$ that are independently distributed between siblings.

Now, under the null hypothesis that $F$ is absent, we note that $P_2$ and $B_1$ are still correlated due to the unblocked path $P_2 - G_p - P_1 - B_1$, so we would expect $\beta \neq 0$ as found. Furthermore, $P_2$ and $B_1$ are still expected to be correlated given $P_1 = 1$ due to the unblocked path $P_2 - G_p - P_1 - E_1 - B_1$, so we would expect $\alpha \neq 0$ as found. Thus, we cannot test the null hypothesis that $F$ is absent without further substantive assumptions beyond those encoded in the causal directed acyclic graph of Figure 1.

Now Hudson et al. [10] were also willing to assume that for no subset of the population did the genetic causes $G_p$ and $G_B$ of $P$ and $B$ prevent disease. Similarly, they assumed there was no subset of the population for whom the environmental causes $E_1$ and $E_2$ of $B$ and $P$ prevented either disease. We will show in Section 5 that under these additional assumptions, the null hypothesis that $F$ is absent implies that the conditional covariance $\alpha$ must be less than or equal to zero, provided that there is no interaction, in the sufficient cause sense, between $E$ and $G_P$. If it is plausible that no sufficient cause interaction between $E$ and $G_P$ exists, then the null hypothesis that $F$ is absent is rejected because the estimate of $\alpha$ is positive with a 95% confidence interval that does not include zero.

Thus, the conclusion in the argument above that familial coaggregation of diseases $B$ and $P$ was present depended critically on the existence of (i) a formal definition of a sufficient cause interaction, (ii) a substantive understanding of what the assumption of no sufficient cause interaction entailed, and (iii) a sound mathematical theory that related assumptions about the absence of sufficient cause interactions to testable restrictions on the distribution of the observed data, specifically on the sign of a particular conditional covariance. In this paper, we provide a theory that offers (i)–(iii).

The remainder of the paper is organized as follows. The second section reviews the directed acyclic graph causal framework and provides some basic definitions; the third section presents the theory which allows for the graphical representation of minimal sufficient causes within the directed acyclic graph causal framework; the fourth section gives an additional preliminary result concerning monotonicity; the fifth section develops results relating minimal sufficient causation and the sign of conditional covariances; the sixth section provides some discussion concerning possible extensions to the present work.

**2. Basic definitions and concepts.** In this section, we review the directed acyclic graph causal framework and give a number of definitions regarding sufficient conjunctions and related concepts. Following Pearl [21], a causal directed acyclic graph is a set of nodes $(X_1, \ldots, X_n)$, corresponding to variables, and directed edges among nodes, such that the graph has no cycles and such that, for each node $X_i$ on the graph, the corresponding variable is given by its nonparametric structural equation $X_i = f_i(pa_i, \epsilon_i)$, where $pa_i$ are the parents of $X_i$ on the graph and the $\epsilon_i$ are mutually independent random variables. These nonparametric structural equations can be seen as a generalization of the path analysis and linear structural equation models [21, 22] developed by Wright [43] in the genetics literature and Haavelmo [9] in the econometrics literature. Robins [27, 28] discusses the close relationship between these nonparametric structural equation models and fully randomized, causally interpreted structured tree graphs [25, 26]. Spirtes, Glymour and Scheines [33] present a causal interpretation of directed acyclic graphs outside the context of nonparametric structural equations and counterfactual variables. It is easily seen from the structural equations that $(X_1, \ldots, X_n)$ admits the following factorization: $p(X_1, \ldots, X_n) = \prod_{i=1}^{n} p(X_i|pa_i)$. The nonparametric structural equations encode counterfactual relationships among the variables represented on the graph. The equations themselves represent one-step ahead counterfactuals with other counterfactuals given by recursive substitution. The requirement that the $\epsilon_i$ be mutually independent is essentially a requirement that there is no variable absent from the graph which, if included on the graph, would be a parent of two or more variables [21, 22].

A path is a sequence of nodes connected by edges regardless of arrowhead direction; a directed path is a path which follows the edges in the direction indicated by the graph's arrows. A node $C$ is said to be a common cause of $A$ and $B$ if there exists a directed path from $C$ to $B$ not through $A$ and a directed path from $C$ to $A$ not through $B$. A collider is a particular node on a path such that both the preceding and subsequent nodes on the path have directed edges going into that node. A backdoor path from $A$ to $B$ is a path that begins with a directed edge going into $A$. A path between $A$ and $B$ is said to be blocked given some set of variables $Z$ if either there is a variable in $Z$ on the path that is not a collider or if there is a collider on the path such that neither the collider itself nor any of its descendants are in $Z$. If all paths between $A$ and $B$ are blocked given $Z$, then $A$ and $B$ are said to be $d$-separated given $Z$. It has been shown that if all paths between $A$ and $B$ are blocked given $Z$, then $A$ and $B$ are conditionally independent given $Z$ [8, 13, 40].

Suppose that a set of nonparametric structural equations represented by a directed acyclic graph $H$ is such that its variables $X$ are partitioned into two sets $X = V \cup W$. If in the nonparametric structural equation for $V \cup W$, by replacing each occurrence of $X_i \in W$ by $f_i(pa_i, \epsilon_i)$, the nonparametric structural equations for $V$ can be written so as to correspond to some causal directed acyclic graph $G$, then $G$ is said to be the marginalization of $H$ over the set of variables $W$. A causal directed acyclic graph with variables $X = V \cup W$ can be marginalized over $W$ if no variable in $W$ is a common cause of any two variables in $V$.

In giving definitions for a sufficient conjunction and related concepts, we will use the following notation. An event is a binary variable taking values in $\{0, 1\}$. The complement of some event $E$ we will denote by $\overline{E}$. A conjunction or product of the events $X_1, \ldots, X_n$ will be written as $X_1 \cdots X_n$. The associative OR operator, $\vee$, is defined by $A \vee B = A + B - AB$. For a random variable $A$ with sample space $\Omega$ we will use the notation $A \equiv 0$ to denote that $A(\omega) = 0$, for all $\omega \in \Omega$. We will use the notation $1_{A=a}$ to denote the indicator function for the random variable $A$ taking the value $a$; for some subset $S$ of the sample space $\Omega$, we will use $1_S$ to denote the indicator that $\omega \in S$. We will use the notation $A \coprod B | C$ to denote that $A$ is conditionally independent of $B$ given $C$. We begin with the definitions of a sufficient conjunction and a minimal sufficient conjunction. These basic definitions make no reference to directed acyclic graphs or causation.

DEFINITION 1. A set of events $X_1, \ldots, X_n$ is said to constitute a sufficient conjunction for event, $D$ if $X_1, \ldots, X_n = 1 \Rightarrow D = 1$.

DEFINITION 2. A set of events $X_1, \ldots, X_n$ which constitutes a sufficient conjunction for $D$ is said to constitute a minimal sufficient conjunction for $D$ if no proper subset of $X_1, \ldots, X_n$ constitutes a sufficient conjunction for $D$.

Sufficient conjunctions for a particular event need not be causes for an event. Suppose a particular sound is produced when and only when an individual blows a whistle. This particular sound the whistle makes is a sufficient conjunction for the whistle's having been blown, but the sound does not cause the blowing of the whistle. The converse, rather, is true; the blowing of the whistle causes the sound to be produced. Corresponding then to these notions of a sufficient conjunction and a minimal sufficient conjunction are those of a sufficient cause and a minimal sufficient cause which will be defined in Section 3.

DEFINITION 3. A set of events $M_1, \ldots, M_n$, each of which may be some product of events, is said to be determinative for some event $D$ if $D = M_1 \vee M_2 \vee \cdots \vee M_n$.

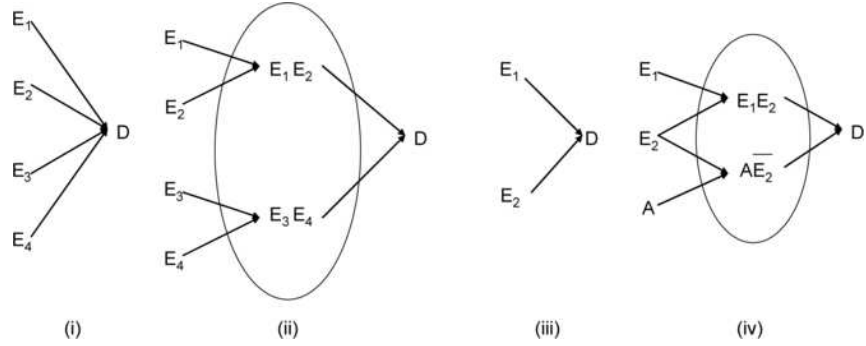

Fig. 2. *Causal directed acyclic graphs with sufficient causation structures.*

DEFINITION 4. A determinative set $M_1, \ldots, M_n$ of (minimal) sufficient conjunctions for $D$ is nonredundant if no proper subset of $M_1, \ldots, M_n$ is determinative for $D$.

EXAMPLE 1. Suppose $A = B \vee CE$ and $D = EF$. If we consider all the minimal sufficient conjunctions for $A$ among the events $\{B, C, D\}$, we can see that $B$ and $CD$ are the only minimal sufficient conjunctions, but it is not the case that $A = B \vee CD$. Clearly then, a complete list of minimal sufficient conjunctions for $A$ generated by a particular collection of events may not be a determinative set of sufficient conjunctions for $A$. If we consider all minimal sufficient conjunctions for $A$ among the events $\{B, C, D, E\}$, we see that $B$ and $CD$ and $CE$ are all minimal sufficient conjunctions. In this example, $B \vee CD \vee CE$ is a determinative set of minimal sufficient conjunctions for $A$ but is not nonredundant. We see then that even when a complete list of minimal sufficient conjunctions generated by a particular collection of events constitutes a determinative set of minimal sufficient conjunctions, it may not be a nonredundant determinative set of minimal sufficient conjunctions.

**3. Minimal sufficient causation and directed acyclic graphs.** In this section, we develop theory which allows for the representation of sufficient conjunctions and minimal sufficient conjunctions on causal directed acyclic graphs. We begin with a motivating example.

EXAMPLE 2. Consider a causal directed acyclic graph given in Figure 2(i). Suppose $E_1 E_2$ and $E_3 E_4$ constitute a determinative set of sufficient conjunctions for $D$. We will show in Theorem 1 below that it follows that the diagram in Figure 2(ii) is also a causal directed acyclic graph where $E_i E_j$ is simply the product or conjunction of $E_i$ and $E_j$; because the sufficient conjunctions $E_1 E_2$ and $E_3 E_4$ are determinative, it follows that $D = E_1 E_2 \vee E_3 E_4$. An ellipse is put around the sufficient conjunctions $E_1 E_2$ and $E_3 E_4$ to indicate that the set is determinative. As will be seen below, in order to add sufficient conjunctions it is important that a determinative set of sufficient conjunctions is known or can be constructed. Consider the causal directed acyclic graph given in Figure 2(iii). Suppose that no determinative set of sufficient conjunctions can be constructed from $E_1$ and $E_2$ alone; suppose further, however, that there exists some other cause of $D$, say $A$, independent of $E_1$ and $E_2$, such that $E_1 E_2$ and $A\overline{E_2}$ form a determinative set of sufficient conjunctions. Then, Theorem 1 below can again be used to show that Figure 2(iv) is a causal directed acyclic graph. Furthermore, it will be shown in Theorem 2 that for any causal directed acyclic graph with a binary node which has only binary parents, a set of variables $\{A_i\}_{i=0}^n$ always exists such that a determinative set of sufficient causes can be formed from the original parents on the graph and the variables $\{A_i\}_{i=0}^n$.

Theorem 1 provides the formal result required for the previous example.

THEOREM 1. *Consider a causal directed acyclic graph $G$ with some node $D$ such that $D$ and all its parents are binary. Suppose that there exists a set of binary variables $A_0, \ldots, A_u$ such that a determinative set of sufficient conjunctions for $D$, say $M_1, \ldots, M_S$, can be formed from conjunctions of $A_0, \ldots, A_u$ along with the parents of $D$ on $G$ and the complements of these variables. Suppose further that there exists a causal directed acyclic graph $H$ such that the parents of $D$ on $H$ that are not on $G$ consist of the nodes $A_0, \ldots, A_u$ and such that $G$ is the marginalization of $H$ over the set of variables which are on the graph for $H$ but not $G$. Then, the directed acyclic graph $J$ formed by adding to $H$ the nodes $M_1, \ldots, M_S$, removing the directed edges into $D$ from the parents of $D$ on $H$, adding directed edges from each $M_i$ into $D$ and adding directed edges into each $M_i$ from every parent of $D$ on $H$ which appears in the conjunction for $M_i$ is itself a causal directed acyclic graph.*

PROOF. To prove that the directed acyclic graph $J$ is a causal directed acyclic graph, it is necessary to show that each of the nodes on the directed acyclic graph can be represented by a nonparametric structural equation involving only the parents on $J$ of that node and a random term $\epsilon_i$ which is independent of all other random terms $\epsilon_j$ in the nonparametric structural equations for the other variables on the graph. The nonparametric structural equation for $M_i$ may be defined as the product of events in the conjunction for $M_i$. The nonparametric structural equation for $D$ can be given by

$$D = M_1 \vee \cdots \vee M_n.$$

The nonparametric structural equations for all other nodes on $J$ can be taken to be the same as those defining the causal directed acyclic graph $H$. Because the nonparametric structural equations for $D$ and for each $M_i$ on $J$ are deterministic, they have no random-error term. Thus, for the nonparametric structural equations defining $D$ and each $M_i$ on $J$, the requirement that the nonparametric structural equation's random term $\epsilon_i$ is independent of all the other random terms $\epsilon_j$ in the nonparametric structural equations for the other variables on the graph is trivially satisfied. That this requirement is satisfied for the nonparametric structural equations for the other variables on $J$ follows from the fact that it is satisfied on $H$. □

In Theorem 1, sufficient conjunctions for $D$ are constructed from some set of variables that, on some causal directed acyclic graph $H$, are all parents of $D$ and thus, within the directed acyclic graph causal framework, it makes sense to speak of sufficient causes and minimal sufficient causes.

DEFINITION 5. If, on a causal directed acyclic graph, some node $D$ with nonparametric structural equation $D = f_D(pa_D, \epsilon_D)$ is such that $D$ and all its parents are binary, then $X_1, \ldots, X_n$ is said to constitute a sufficient cause for $D$ if $X_1, \ldots, X_n$ are all parents of $D$ or complements of the parents of $D$ and are such that $f_D(pa_D, \epsilon_D) = 1$ for all $\epsilon_D$ whenever $pa_D$ is such that $X_1 \cdots X_n = 1$; if no proper subset of $X_1, \ldots, X_n$ also constitutes a sufficient cause for $D$, then $X_1, \ldots, X_n$ is said to constitute a minimal-sufficient cause for $D$. A set of (minimal) sufficient causes, $M_1, \ldots, M_n$, each of which is a product of the parents of $D$ and their complements, is said to be determinative for some event $D$ if, for all $\epsilon_D$, $f_D(pa_D, \epsilon_D) = 1$ if and only if $pa_D$ is such that $M_1 \vee M_2 \vee \cdots \vee M_n = 1$; if no proper subset of $M_1, \ldots, M_n$ is also determinative for $D$, then $M_1, \ldots, M_n$ is said to constitute a nonredundant determinative set of (minimal) sufficient causes for $D$.

If, for some directed acyclic graph $G$ there exist $A_0, \ldots, A_u$ which satisfy the conditions of Theorem 1 for some node $D$ on $G$ so that a determinative set of sufficient causes for $D$ can be constructed from $A_0, \ldots, A_u$ along with the parents of $D$ on $G$ and their complements, then $D$ will be said to admit a sufficient causation structure. As in Example 2, we will, in general, replace the $M_i$ nodes with the conjunctions that constitute them. The node $D$ with directed edges from the $M_i$ nodes is effectively an OR node. The $M_i$ nodes with the directed edges from the $A_i$ nodes and the parents of $D$ on $G$ are effectively AND nodes. We call this resulting diagram a causal directed acyclic graph with a sufficient causation structure (or a minimal sufficient causation structure if the determinative set of sufficient conjunctions for $D$ are each minimal sufficient conjunctions).

Because a causal directed acyclic graph with a sufficient causation structure is itself a causal directed acyclic graph, the $d$-separation criterion applies and allows one to determine independencies and conditional independencies. A minimal sufficient causation structure will often make apparent conditional independencies within a particular stratum of the conditioning variable which were not apparent on the original causal directed acyclic graph. The following corollary is useful in this regard.

COROLLARY 1. *If some node $D$ on a causal directed acyclic graph admits a sufficient causation structure then conditioning on $D = 0$ conditions also all sufficient cause nodes for $D$ on the causal directed acyclic graph with the sufficient causation structure.*

EXAMPLE 2 (Continued). Consider the causal directed acyclic graph with the minimal sufficient causation structure given in Figure 2(ii). Conditioning on $D = 0$ also conditions on $E_1 E_2 = 0$ and $E_3 E_4 = 0$, and thus, by the $d$-separation criterion, $E_i$ is conditionally independent of $E_j$ given $D = 0$ for $i \in \{1, 2\}, j \in \{3, 4\}$. In the causal directed acyclic graph with the minimal sufficient causation structure in Figure 2(iv), no similar conditional independence relations within the $D = 0$ stratum holds. Although conditioning on $D = 0$ conditions also on $E_1 E_2 = 0$ and $A\overline{E_2} = 0$ there still remains an unblocked path $E_1 - E_1 E_2 - E_2 - A\overline{E_2} - A$ between $E_1$ and $A$, and so $E_1$ and $A$ are not conditionally independent given $D = 0$; Similarly, there are unblocked paths between $E_1$ and $E_2$ given $D = 0$ and also between $E_2$ and $A$ given $D = 0$.

The additional variables $A_0, \ldots, A_u$ needed to form a set of sufficient causes for $D$ we will refer to as the co-causes of $D$. The co-causes $A_0, \ldots, A_u$ required to form a determinative set of sufficient conjunctions for $D$ will generally not be unique. For example, if $D = A_0 \vee A_1 E$ then it is also the case that $D = B_0 \vee B_1 E$, where $B_0 = A_0$ and $B_1 = \overline{A_0} A_1$. Similarly, there will, in general, be no unique set of sufficient causes that is determinative for $D$. For example, if $E_1$ and $E_2$ constitute a set of sufficient causes for $D$ so that $D = E_1 \vee E_2$, then it is also the case that $E_1 \overline{E_2}$, $\overline{E_1} E_2$, and $E_1 E_2$ also constitute a set of sufficient causes for $D$, and so we could also write $D = E_1 \overline{E_2} \vee \overline{E_1} E_2 \vee E_1 E_2$. It can be shown that not even nonredundant determinative sets of minimal sufficient causes are unique.

Corresponding to the definition of a sufficient cause is the more philosophical notion of a causal mechanism. A causal mechanism can be conceived of as a set of events or conditions which, if all present, bring about the outcome under consideration through a particular pathway. A causal mechanism thus provides a particular description of how the outcome comes about. Suppose, for instance, that an individual were exposed to two poisons, $E_1$ and $E_2$, such that in the absence of $E_2$, the poison $E_1$ would lead to heart failure resulting in death; and that in the absence of $E_1$, the poison $E_2$ would lead to respiratory failure resulting in death; but such that when $E_1$ and $E_2$ are both present, they interact and lead to a failure of the nervous system again resulting in death. In this case, there are three distinct causal mechanisms for death each corresponding to a sufficient cause for $D$: death by heart failure corresponding to $E_1\overline{E_2}$, death by respiratory failure corresponding to $\overline{E_1}E_2$ and death due to a failure of the nervous system corresponding to $E_1E_2$. It is interesting to note that in this case none of the sufficient causes corresponding to the causal mechanisms is minimally sufficient. Each of $E_1\overline{E_2}$, $\overline{E_1}E_2$ and $E_1E_2$ is sufficient for $D$ but none is minimally sufficient, as either $E_1$ or $E_2$ alone is sufficient for death. We will refer to a sufficient cause for $D$ as a causal mechanism for $D$ if the node for the sufficient cause corresponds to a variable, potentially subject to intervention, which whenever the variable takes the value 1, the outcome $D$ inevitably results.

The last example shows that the existence of a particular set of determinative sufficient causes does not guarantee that there are actual causal mechanisms corresponding to these sufficient causes; it only implies that a set of causal mechanisms corresponding to these sufficient causes cannot be ruled out by a complete knowledge of counterfactual outcomes. In particular, in the previous example, the set $\{E_1, E_2\}$ is a determinative set of sufficient causes that does not correspond to the actual set of causal mechanisms $\{E_1\overline{E_2}, \overline{E_1}E_2, E_1E_2\}$. If there are two or more sets of sufficient causes that are determinative for some outcome $D$ then although the two sets of determinative sufficient causes are logically equivalent for prediction, we nevertheless view them as distinct. In such cases, some knowledge of the subject matter in question will, in general, be needed to discern which of the sets of determinative sufficient causes actually corresponds to the true causal mechanisms. For instance, in the previous example, we needed biological knowledge of how poisons brought about death in the various scenarios. We will, in the interpretation of our results, assume that there always exists some set of true causal mechanisms which forms a determinative set of sufficient causes for the outcome. The concept of synergism is closely related to that of a causal mechanism and is often found in the epidemiologic literature [11, 29, 32]. We will say that there is synergism between the effects of $E_1$ and $E_2$ on $D$ if there exists a sufficient cause for $D$ which represents some causal mechanism and such that this sufficient cause has $E_1$ and $E_2$ in its conjunction. In related work, we have developed tests for synergism, that is, tests for the joint presence of two or more causes in a single sufficient cause [36, 37]. In some of our examples and in our discussion of the various results in the paper, we will sometimes make reference to the concepts of a causal mechanism and synergism. However, all definitions, propositions, lemmas, theorems and corollaries will be given in terms of sufficient causes for which we have a precise definition.

The graphical representation of sufficient causes on a causal directed acyclic graph does not require that the determinative set of sufficient causes for $D$ be minimally sufficient, nor does it require that the set of determinative sufficient causes for $D$ be nonredundant. To expand a directed acyclic graph into another directed acyclic graph with sufficient cause nodes, all that is required is that the set of sufficient causes constitutes a determinative set of sufficient causes for $D$. However, a set of events that constitutes a sufficient cause can be reduced to a set of events that constitutes a minimal sufficient cause by iteratively excluding unnecessary events from the set until a minimal sufficient cause is obtained. Also, a set of determinative sufficient causes that is redundant can be reduced to one that is nonredundant by excluding those sufficient causes or minimal sufficient causes that are redundant. It is sometimes an advantage to reduce a redundant set of sufficient causes to a nonredundant set of minimal sufficient causes. This is so because allowing sufficient causes that are not minimally sufficient or allowing redundant sufficient causes or redundant minimal sufficient causes can obscure the conditional independence relations implied by the structure of the causal directed acyclic graph. This is made evident in Example 3.

EXAMPLE 3. Consider the causal directed acyclic graph with the minimal sufficient causation structure given in Figure 3(i). Conditioning on $D = 0$ conditions also on $AB = 0$ and $EF = 0$ and by the $d$-separation criterion, $A$ and $E$ are conditionally independent given $D = 0$. But now consider an expanded structure for this causal directed acyclic graph which involves only minimal sufficient causes but which allows redundant minimal sufficient causes. Define $Q = BE$, then $AQ$ is a minimal sufficient cause for $D$ since $AQ = 1 \Rightarrow AB = 1 \Rightarrow D = 1$, but $A = 1 \nRightarrow D = 1$ and $Q = 1 \nRightarrow D = 1$. Now $AB, AQ, EF$ is a determinative but redundant set of minimal sufficient causes for $D$. Figure 3(ii) gives an alternative causal directed acyclic graph with a minimal sufficient causation structure for the causal relationships indicated in Figure 3(i). In Figure 3(ii), conditioning on $D = 0$ conditions also on $AB = 0$, $AQ = 0$ and $EF = 0$, but the $d$-separation criteria no longer imply that $A$ and $E$ are conditionally independent given $D = 0$; because of conditioning on $D = 0$, there is an unblocked path between $A$ and $E$, namely $A - AQ - Q - BE - E$. Allowing the redundant minimal sufficient cause $AQ$ in the minimal sufficient causation structure obscures the conditional independence relation. Similar examples can be constructed to show that allowing sufficient causes that are not minimally sufficient can also obscure conditional independence relations [35].

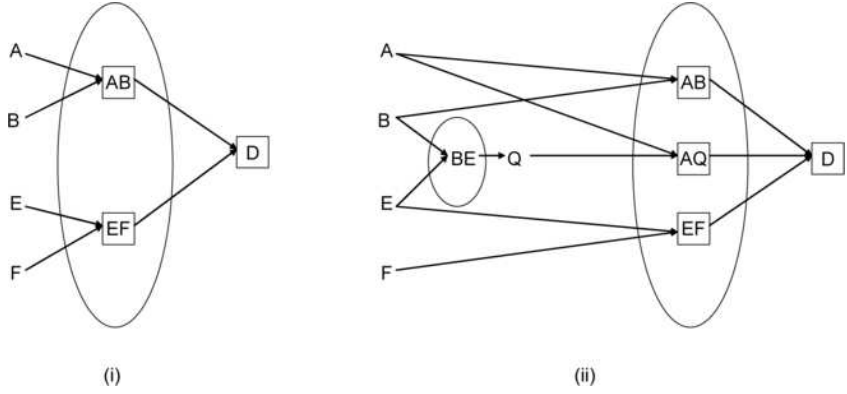

Fig. 3. *Example illustrating that redundant sufficient causes can obscure conditional independence relations.*

Although allowing sufficient causes that are not minimally sufficient or allowing redundant sufficient causes or redundant minimal sufficient causes can obscure the conditional independence relations implied by the structure of the causal directed acyclic graph, it may sometimes be desirable to include nonminimal sufficient causes or redundant sufficient causes. For example, as noted above, nonminimal sufficient cause nodes or redundant sufficient cause nodes may represent separate causal mechanisms upon which it might be possible to intervene. Further discussion of conditional independence relations in sufficient causation structures with nonminimally sufficient causes and redundant sufficient causes is given in Section 6.

Note a sufficient cause need only involve one co-cause $A_i$ in its conjunction because if it involved $A_{i_1}, \ldots, A_{i_k}$, then $A_{i_1}, \ldots, A_{i_k}$ could be replaced by the product $A'_i = A_{i_1} \cdots A_{i_k}$. In certain cases though, it may be desirable to include more than one $A_i$ in a sufficient cause if this corresponds to the actual causal mechanisms. If a set of variables $A_0, \ldots, A_u$ satisfying Theorem 1 can be constructed from functions of the random term $U = \epsilon_D^G$ of the nonparametric structural equation for $D$ on $G$ and their complements so that $A_i = f_i(U)$, then $H$ can be chosen to be the graph $G$ with the additional nodes $U, A_0, \ldots, A_u$ and with directed edges from $U$ into each $A_i$ and from each $A_i$ into $D$. This gives rise to the definition, given below, of a representation for $D$.

DEFINITION 6. If $D$ and all of its parents on the causal directed acyclic graph $G$ are binary and there exists some set $\{A_i, P_i\}$ such that each $P_i$ is some conjunction of the parents of $D$ and their complements, such that there exist functions $f_i$ for which $A_i = f_i(\epsilon_D)$, where $\epsilon_D$ is the random term in the nonparametric structural equation for $D$ on $G$ and such that $D = \bigvee_i A_i P_i$, then $\{A_i, P_i\}$ is said to constitute a representation for $D$.

If the $A_i$ variables are constructed from functions of the random term $\epsilon_D$ in the nonparametric structural equation for $D$ on $G$, then these $A_i$ variables may or may not allow for interpretation, and they may or may not be such that an intervention on these $A_i$ variables is conceivable. In certain cases, the $A_i$ variables may simply be logical constructs for which no intervention is conceivable. Although in certain cases it may not be possible to intervene on the $A_i$ variables, we will still refer to conjunctions of the form $A_i P_i$ as sufficient *causes* for $D$, as it is assumed that it is possible to intervene on the parents of $D$ which constitute the conjunction for $P_i$.

Suppose that for some node $D$ on a causal directed acyclic graph $G$, a set of variables $A_0, \ldots, A_u$ satisfying Theorem 1 can be constructed from functions of the random term $U = \epsilon_D$ in the nonparametric structural equation for $D$ on $G$, so that a representation for $D$ is given by $D = \bigvee_i A_i P_i$. Then, in order to simplify the diagram, instead of adding to $G$ the variable $U$ and directed edges from $U$ into each $A_i$ so as to form the minimal sufficient causation structure, we will sometimes suppress $U$ and simply add an asterisk next to each $A_i$ indicating that the $A_i$ variables have a common cause.

PROPOSITION 1. *For any representation for $D$, the co-causes $A_i$ will be independent of the parents of $D$ on the original directed acyclic graph $G$.*

PROOF. This follows immediately from the fact that for any representation for $D$, the co-causes are functions of the random term in the nonparametric structural equation for $D$. □

If some of the sufficient causes for $D$ are unknown, then it is not obvious how one might make use of Theorem 1. The theorem allowed for a sufficient causation structure on a causal directed acyclic graph, provided there existed some set of co-causes $A_0, \ldots, A_u$. Theorem 2 complements Theorem 1 in that it essentially states that when $D$ and all of its parents are binary such a set of co-causes always exists. The variables $A_0, \ldots, A_u$ are constructed from functions of the random term $\epsilon_D$ in the nonparametric structural equation for $D$ on $G$. Before stating and proving Theorem 1, we illustrate how the co-causes can be constructed by a simple example.

EXAMPLE 4. Suppose $E$ is the only parent of $D$, then the structural equation for $D$ is given by $D = f(E, \varepsilon_D)$. Define $A_0$, $A_1$ and $A_2$ as follows: let $A_0(\omega) = 1$ if $f(1, \varepsilon_D(\omega)) = f(0, \varepsilon_D(\omega)) = 1$ and $A_0(\omega) = 0$ otherwise; let $A_1(\omega) = 1$ if $f(1, \varepsilon_D(\omega)) = 1$ and $f(0, \varepsilon_D(\omega)) = 0$, and $A_1(\omega) = 0$ otherwise; and let $A_2(\omega) = 1$ if $f(1, \varepsilon_D(\omega)) = 0$ and $f(0, \varepsilon_D(\omega)) = 1$, and $A_2(\omega) = 0$ otherwise. It is easily verified that $D = A_0 \vee A_1 E \vee A_2 \overline{E}$ and that $A_0$, $A_1 E$ and $A_2 \overline{E}$ constitute a determinative set of minimal sufficient causes for $D$. Note that this construction will give a determinative set of minimal sufficient causes for $D$ regardless of the form of $f$ and the distribution of $\varepsilon_D$.

MINIMAL SUFFICIENT CAUSATION 15THEOREM 2. *Consider a causal directed acyclic graph $G$ on which there exists some node $D$ such that $D$ and all its parents are binary, then there exist variables $A_0, \ldots, A_u$ that satisfy the conditions of Theorem 1 and such that the sufficient causes constructed from $A_0, \ldots, A_u$ along with the parents of $D$ on $G$ and their complements are, in fact, minimal sufficient causes.*

PROOF. The nonparametric structural equation for $D$ is given by $D = f(pa_D, \varepsilon_D)$. Suppose $D$ has $m$ parents on the original causal directed acyclic graph $G$. Since these parents are binary, there are $2^m$ values which $pa_D$ can take. Since $f$ maps $(pa_D, \varepsilon_D)$ to $\{0, 1\}$, each value of $\varepsilon_D$ assigns to every possible realization of $pa_D$ either 0 or 1 through $f$. There are $2^{2^m}$ such assignments. Thus, without loss of generality, we may assume that $\varepsilon_D$ takes on some finite number of distinct values $N \leq 2^{2^m}$; and so, we may write the sample space for $\varepsilon_D$ as $\Omega_D = \{\omega_1, \ldots, \omega_N\}$, and we may use $\omega = \omega_i$ and $\varepsilon_D = \varepsilon_D(\omega_i)$ interchangeably. The co-causes $A_0, \ldots, A_u$ can be constructed as follows. Let $W_i$ be the indicator $1_{\varepsilon_D = \varepsilon_D(\omega_i)}$. Let $P_i$ be some conjunction of the parents of $D$ and their complements, that is, $P_i = F_1^i \cdots F_{n_i}^i$, where each $F_k^i$ is either a parent of $D$, say $E_j$ or its complement $\overline{E_j}$. For each $P_i$, let $A_i \equiv 1$ if $F_1^i \cdots F_{n_i}^i$ is a minimal sufficient cause for $D$ and

$$A_i = \bigvee_j \{W_j : W_j F_1^i \cdots F_{n_i}^i \text{ is a minimal sufficient cause for } D\}$$

otherwise. Let $M_i = P_i$ if $A_i = 1$, and $M_i = A_i P_i$ otherwise. It must be shown that each $M_i = A_i F_1^i \cdots F_{n_i}^i$ is a minimal sufficient cause and that the set of $M_i$'s constitutes a minimal sufficient cause representation for $D$ (or more precisely, the set of $M_i$'s for which $A_i$ is not identically 0 constitutes a minimal sufficient cause representation for $D$). We first show that each $M_i = A_i F_1^i \cdots F_{n_i}^i$ is a minimal sufficient cause for $D$. Clearly, this is the case if $A_i \equiv 1$. Now consider those $A_i$ such that $A_i$ is not identically 0 and not identically 1 and suppose $A_i = W_1^i \vee \cdots \vee W_{v_i}^i$, where each $W_j^i$ is such that $W_j^i F_1^i \cdots F_{n_i}^i$ is a minimal sufficient cause for $D$. If $A_i F_1^i \cdots F_{n_i}^i$ is not a minimal sufficient cause, then either $F_1^i \cdots F_{n_i}^i = 1 \Rightarrow D = 1$ or there exists $j$ such that

$$A_i F_1^i \cdots F_{j-1}^i F_{j+1}^i \cdots F_{n_i}^i \quad \Rightarrow \quad D = 1.$$

Suppose first that $F_1^i \cdots F_{n_i}^i = 1 \Rightarrow D = 1$ then there does not exist a $W_j$ such that $W_j F_1^i \cdots F_{n_i}^i$ is a minimal sufficient cause for $D$; but this contradicts $A_i$ is not identically 1. On the other hand, if there exists $j$ such that $A_i F_1^i \cdots F_{j-1}^i F_{j+1}^i \cdots F_{n_i}^i \Rightarrow D = 1$, then it is also the case that

$$W_1^i F_1^i \cdots F_{j-1}^i F_{j+1}^i \cdots F_{n_i}^i \quad \Rightarrow \quad D = 1,$$

since $A_i$ is simply a disjunction of the $W_j^i$'s. However, it would then follow that $W_1^i F_1^i \cdots F_{n_i}^i$ is not a minimal sufficient cause for $D$. But this contradicts the definition of $W_1^i$. Thus, $A_i F_1^i \cdots F_{n_i}^i$ must be a minimal sufficient cause for $D$. It remains to be shown that the set of $M_i$'s for which $A_i$ is not identically 0 constitutes a minimal sufficient cause representation for $D$. We must show that if $D = 1$, then there exists a $M_i = A_i P_i$ for which $M_i = 1$. Now $D$ is a function of $(\varepsilon_D, E_1, \ldots, E_m)$, so let $(\varepsilon_D^*, E_1^*, \ldots, E_m^*)$ be any particular value of $(\varepsilon_D, E_1, \ldots, E_m)$ for which $D = 1$. Consider the set $\{E_1, \ldots, E_m\}$. If for any $j$,

$$\varepsilon_D = \varepsilon_D^*, \qquad E_1 = E_1^*, \ldots, E_{j-1} = E_{j-1}^*,$$
$$E_{j+1} = E_{j+1}^*, \ldots, E_m = E_m^* \quad \Rightarrow \quad D = 1,$$

remove $E_j$ from $\{E_1, \ldots, E_m\}$. Continue to remove those $E_j$ from this set which are not needed to maintain the implication $D = 1$. Suppose the set that remains is $\{E_{h_1}, \ldots, E_{h_S}\}$, then either we have $E_{h_1} = E_{h_1}^*, \ldots, E_{h_S} = E_{h_S}^* \Rightarrow D = 1$ or we have

$$E_{h_1} = E_{h_1}^*, \ldots, E_{h_S} = E_{h_S}^* \not\Rightarrow D = 1$$

and

$$\varepsilon_D = \varepsilon_D^*, \qquad E_{h_1} = E_{h_1}^*, \ldots, E_{h_S} = E_{h_S}^* \quad \Rightarrow \quad D = 1.$$

If $E_{h_1} = E_{h_1}^*, \ldots, E_{h_S} = E_{h_S}^* \Rightarrow D = 1$, then if we define $F_j$ as the indicator $F_j = 1_{(E_{h_j} = E_{h_j}^*)}$, $F_1 \cdots F_S$ is a minimal sufficient cause for $D$ and there thus exists an $i$, such that $P_i = F_1 \cdots F_S$ and $M_i = P_i$, and when $E_{h_1} = E_{h_1}^*, \ldots, E_{h_S} = E_{h_S}^*$, we have $M_i = 1$. If $E_{h_1} = E_{h_1}^*, \ldots, E_{h_S} = E_{h_S}^* \not\Rightarrow D = 1$ but $\varepsilon_D = \varepsilon_D^*, E_{h_1} = E_{h_1}^*, \ldots, E_{h_S} = E_{h_S}^* \Rightarrow D = 1$, then if we define $F_j$ as the indicator $1_{(E_{h_j} = E_{h_j}^*)}$, $1_{\varepsilon_D = \varepsilon_D^*} F_1 \cdots F_S$ is a minimal sufficient cause for $D$; and there exists an $i$ such that $M_i = A_i P_i$ and $P_i = F_1 \cdots F_S$; and $\varepsilon_D = \varepsilon_D^* \Rightarrow A_i = 1$, such that

$$\varepsilon_D = \varepsilon_D^*, \qquad E_{h_1} = E_{h_1}^*, \ldots, E_{h_S} = E_{h_S}^* \quad \Rightarrow \quad M_i = 1.$$

We have thus shown when $D = 1$, there exists an $M_i$ such that $M_i = 1$ and so the $M_i$'s constitutes a minimal sufficient cause representation for $D$. □

The variables $A_i$ constructed in Theorem 2, along with their corresponding conjunctions $P_i$ of the parents of $D$ and their complements, we define below as the canonical representation for $D$. It is easily verified that the co-causes and representation constructed in Example 4 is the canonical representation for $D$ in that example.

DEFINITION 7. Consider a causal directed acyclic graph $G$, such that some node $D$ and all of its parents are binary. Let $\Omega_D$ be the sample space for the random term $\epsilon_D$ in the nonparametric structural equation for $D$ on $G$. The conjunctions $P_i = F_1^i \cdots F_{n_i}^i$, where each $F_k^i$ is either a parent of $D$ or the complement of a parent of $D$, along with the variables $A_i$ constructed by $A_i \equiv 1$ if $F_1^i \cdots F_{n_i}^i$ is a minimal sufficient cause for $D$ and $A_i = \bigvee_{\omega_j \in \Omega_D} \{1_{\varepsilon_D = \varepsilon_D(\omega_j)} : 1_{\varepsilon_D = \varepsilon_D(\omega_j)} F_1^i \cdots F_{n_i}^i$ is a minimal sufficient cause for $D\}$; otherwise, is said to be the canonical representation for $D$.

As noted above, there will in general exist more than one set of co-causes $A_0, \ldots, A_u$, which together with the parents of $D$ and their complements can be used to construct a sufficient cause representation for $D$. The set of $A_i$'s in the canonical representation constitutes only one particular set of variables which can be used to construct a sufficient cause representation. If $D$ has three or more parents, examples can be constructed in which the canonical representation is redundant. Examples can also be constructed to show that when the canonical representation is redundant, it is not always *uniquely* reducible to a nonredundant minimal sufficient cause representation. Although the canonical representation will not always be nonredundant, it does however guarantee that for a binary variable with binary parents, a determinative set of minimal sufficient causes always exists. The canonical representation in a sense "favors" conjunctions with fewer terms. As can be seen in the simple illustration given in Example 4, the canonical representation will never have $A_i = 1$, for some conjunction $P_i$, when there is a conjunction $P_j$ with $A_j = 1$ and such that the components of $P_j$ are a subset of those in the conjunction for $P_i$.

**4. Monotonic effects and minimal sufficient causation.** Minimal sufficient causes for a particular event $D$ may have present in their conjunction the parents of $D$ or the complements of these parents. In certain cases, no minimal sufficient cause will involve the complement of a particular parent of $D$. Such cases closely correspond to what will be defined below as a positive monotonic effect. Essentially, a positive monotonic effect will be said to be present when a function in a nonparametric structural equation is nondecreasing in a particular argument for all values of the other arguments of the function. In this section, we develop the relationship between minimal sufficient causation and monotonic effects.

DEFINITION 8. The nonparametric structural equation for some node $D$ on a causal directed acyclic graph with parent $E$ can be expressed as $D = f(\widetilde{pa}_D, E, \epsilon_D)$, where $\widetilde{pa}_D$ are the parents of $D$ other than $E$; $E$ is said to have a positive monotonic effect on $D$ if, for all $\widetilde{pa}_D$ and $\epsilon_D$, $f(\widetilde{pa}_D, E_1, \epsilon_D) \geq$ $f(\widetilde{pa}_D, E_2, \epsilon_D)$ whenever $E_1 \geq E_2$. Similarly, $E$ is said to have a negative monotonic effect on $D$ if, for all $\widetilde{pa}_D$ and $\epsilon_D$, $f(\widetilde{pa}_D, E_1, \epsilon_D) \leq f(\widetilde{pa}_D, E_2, \epsilon_D)$ whenever $E_1 \geq E_2$.

Note that this notion of a monotonic effect is somewhat stronger than Wellman's qualitative probabilistic influence [41]. See [38, 39] for further discussion.

THEOREM 3. *If $E$ is parent of $D$ and if $D$ and all of its parents are binary, then the following are equivalent:* (i) *$E$ has a positive monotonic effect on $D$;* (ii) *there is some representation for $D$ which is such that none of the representation's conjunctions contain $\overline{E}$;* (iii) *the canonical representation of $D$, $\bigvee_i A_i P_i$, is such that no conjunction $P_i$ contains $\overline{E}$.*

PROOF. We see that (iii) implies (ii) because the representation required by (ii) is met by the canonical representation of $D$, as constructed in Theorem 2. To show that (ii) implies (i), we assume that we have a representation for $D$ such that $D = \bigvee_i A_i P_i$, where each $P_i$ is some conjunction of the parents of $D$ and their complements but does not contain $\overline{E}$. If $f(\widetilde{pa}_D, \overline{E}, \epsilon_D) = 1$, then $f(\widetilde{pa}_D, E, \epsilon_D) = 1$ because $D = \bigvee_i A_i P_i$ and none of the $P_i$ involve $\overline{E}$; from this, (i) follows. To show that (i) implies (iii) we prove the contrapositive. Suppose that the canonical representation of $D$, $\{A_i, P_i\}$, is such that there exists a $P_i$ which contains $\overline{E}$ in its conjunction. Then there exists some value $\varepsilon_D^*$ of $\varepsilon_D$ and some conjunction of the parents of $D$ and their complements, say $F_1 \cdots F_n$, such that $W_i F_1 \cdots F_n \overline{E}$ constitutes a minimal sufficient cause for $D$, where $W_i = 1_{(\varepsilon_D^* = \varepsilon_D)}$. Let $\widetilde{pa}_D^*$ take the values given by $F_1 \cdots F_n$. This may not suffice to fix $\widetilde{pa}_D^*$, but there must exist some value of the remaining parents of $D$ other than $E$ which, in conjunction with $W_i F_1 \cdots F_n E$, gives $D = 0$; for if there were no such values of the other parents, then $W_i F_1 \cdots F_n$ itself would be sufficient for $D$, and $W_i F_1 \cdots F_n \overline{E}$ would not be a minimal sufficient cause for $D$. Let $\widetilde{pa}_D^*$ be such that $\widetilde{pa}_D^*$ and $\overline{E}$ together with $\varepsilon_D^*$ give $D = 1$, but $\widetilde{pa}_D^*$ and $E$ with $\varepsilon_D^*$ give $D = 0$. Then, $f(\widetilde{pa}_D^*, \overline{E}, \varepsilon_D^*) = 1$, but $f(\widetilde{pa}_D^*, E, \varepsilon_D^*) = 0$, and thus, (i) does not hold. This completes the proof. □

**5. Conditional covariance and minimal sufficient causation.** When two binary parents of some event $D$ have positive monotonic effects on $D$, it is in some cases possible to determine the sign of the conditional covariance of these two parents. In general, even in the setting of monotonic effects, the conditional covariance may be of either positive or negative sign; however, when additional knowledge is available concerning the minimal sufficient causation structure of $D$, it is often possible to determine the sign of the conditional covariance of two parents of $D$. Theorem 4 gives conditions under which the sign of the conditional covariance can be determined. Theorems 5 and 6 extend the conclusions of Theorem 4 to certain cases concerning the conditional covariance of two variables that may not be parents of the conditioning variable. The proof of Theorem 4 is suppressed; the proof involves extensive but routine algebraic manipulation and factoring (details are available from the authors upon request).

THEOREM 4. *Suppose that $E_1$ and $E_2$ are the only parents of $D$ on some causal directed acyclic graph, that $E_1$, $E_2$ and $D$ are all binary and that both $E_1$ and $E_2$ have a positive monotonic effect on $D$. Then, for any representation for $D$ such that $D = A_0 \vee A_1 E_1 \vee A_2 E_2 \vee A_3 E_1 E_2$, the following hold:*

(i) *If $A_0 \equiv 0$, then $\mathrm{Cov}(E_1, E_2|D) \leq 0$.*
(ii) *If $A_0 \equiv 0$, $A_1$ and $A_2$ are independent and $E_1$ and $E_2$ are independent, then $\mathrm{Cov}(E_1, E_2|\overline{D}) \leq 0$.*
(iii) *If $A_1 \equiv 1$ or $A_2 \equiv 1$, then $\mathrm{Cov}(E_1, E_2|D) \leq 0$ provided $\mathrm{Cov}(E_1, E_2) \leq 0$.*
(iv) *If $A_1 \equiv 1$ or $A_2 \equiv 1$, then $\mathrm{Cov}(E_1, E_2|\overline{D}) = 0$.*
(v) *If $A_1 \equiv 0$ or $A_2 \equiv 0$, then $\mathrm{Cov}(E_1, E_2|D) \geq 0$ provided $\mathrm{Cov}(E_1, E_2) \geq 0$.*
(vi) *If $A_1 \equiv 0$ or $A_2 \equiv 0$, then $\mathrm{Cov}(E_1, E_2|\overline{D}) \leq 0$ provided $\mathrm{Cov}(E_1, E_2) \leq 0$.*
(vii) *If $A_3 \equiv 0$, then $\mathrm{Cov}(E_1, E_2|D) \leq 0$ provided $\mathrm{Cov}(E_1, E_2) \leq 0$.*
(viii) *If $A_3 \equiv 0$, $A_1$ and $A_2$ are independent, $E_1$ and $E_2$ are independent and also $A_0$ is independent of either $A_1$ or $A_2$, then $\mathrm{Cov}(E_1, E_2|\overline{D}) = 0$.*

Note that parts (i)–(viii) of Theorem 4 all require some knowledge of a sufficient cause representation for $D$, that is, that $A_0 = 0$, $A_1 \equiv 1$ or $A_1 \equiv 0$, etc. Conclusions about the sign of the conditional covariance cannot be drawn from Theorem 4 without some knowledge of a sufficient causation structure. In general, this knowledge of a sufficient causation structure would come from prior beliefs about the actual causal mechanisms for $D$. As can be seen from Theorem 4, if no knowledge of the sufficient causes is available, the conditional covariances $\mathrm{Cov}(E_1, E_2|D)$ and $\mathrm{Cov}(E_1, E_2|\overline{D})$ may be of either sign, even if $E_1$ and $E_2$ have positive monotonic effects on $D$. For example, if $E_1$ and $E_2$ have positive monotonic effects on $D$ and (v) holds then $\mathrm{Cov}(E_1, E_2|D) \geq 0$; but if $E_1$ and $E_2$ have positive monotonic effects on $D$ and (i) holds, then $\mathrm{Cov}(E_1, E_2|D) \leq 0$.

If $E_1$ and $E_2$ are the only parents of $D$, possibly correlated due to some common cause $C$, and have positive monotonic effects on $D$ then the minimal sufficient causation structure for the causal directed acyclic graph is that given in Figure 4.

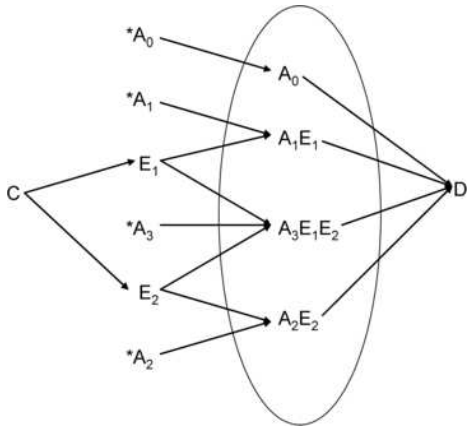

Fig. 4. *Minimal sufficient causation structure when $E_1$ and $E_2$ have positive monotonic effects on $D$.*

Recall the asterisk is used to indicate that $A_0$, $A_1$, $A_2$ or $A_3$ may have a common cause $U$. If one of $A_0$, $A_1$, $A_2$ or $A_3$ is identically 0 or 1, then Theorem 4 may be used to draw conclusions about the sign of the conditional covariance $\text{Cov}(E_1, E_2|D)$. For example, if one believes that there is no synergism between $E_1$ and $E_2$ in the actual causal mechanisms for $D$ then $A_3 \equiv 0$; if this holds, then parts (vii) and (viii) of Theorem 4 can be used to determine the sign of the conditional covariance. Theorem 4 has an obvious analogue if one or both of $E_1$ or $E_2$ have a negative monotonic effect on $D$. If $D$ has more than two parents, but if the two parents, $E_1$ and $E_2$, are independent of all other parents of $D$, then the causal directed acyclic graph can be marginalized over these other parents, and Theorem 4 could be applied to the resulting causal directed acyclic subgraph.

Some of the conclusions of Theorem 4 require knowing the sign of $\text{Cov}(E_1, E_2)$ and Proposition 2 below (proved elsewhere [39]) relates the sign of $\text{Cov}(E_1, E_2)$ to the presence of monotonic effects. In order to state this proposition and to allow for the development of extensions to Theorem 4, we need a few additional definitions.

DEFINITION 9. An edge on a causal directed acyclic graph from $X$ to $Y$ is said to be of positive (negative) sign if $X$ has a positive (negative) monotonic effect on $Y$. If $X$ has neither a positive monotonic effect nor a negative monotonic effect on $Y$, then the edge from $X$ to $Y$ is said to be without a sign.

DEFINITION 10. The sign of a path on a causal directed acyclic graph is the product of the signs of the edges that constitute that path. If one of the edges on a path is without a sign, then the sign of the path is said to be undefined.

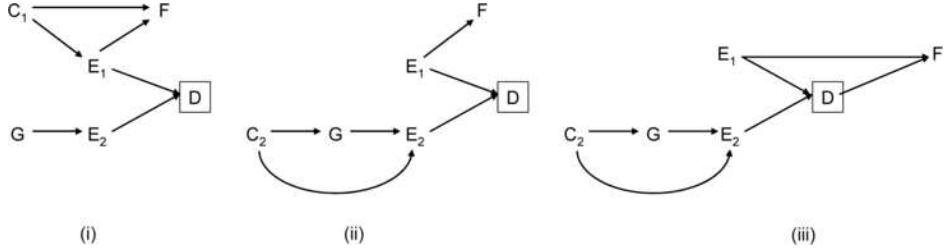

Fig. 5. *Examples requiring extensions to Theorem 4.*

DEFINITION 11. Two variables $X$ and $Y$ are said to be positively monotonically associated if all directed paths between $X$ and $Y$ are of positive sign, and all common causes $C_i$ of $X$ and $Y$ are such that all directed paths from $C_i$ to $X$ not through $Y$ are of the same sign as all directed paths from $C_i$ to $Y$ not through $X$; the variables $X$ and $Y$ are said to be negatively monotonically associated if all directed paths between $X$ and $Y$ are of negative sign, and all common causes $C_i$ of $X$ and $Y$ are such that all directed paths from $C_i$ to $X$ not through $Y$ are of the opposite sign as all directed paths from $C_i$ to $Y$ not through $X$.

PROPOSITION 2. *If $X$ and $Y$ are positively monotonically associated, then $\mathrm{Cov}(X,Y) \geq 0$. If $X$ and $Y$ are negatively monotonically associated, then $\mathrm{Cov}(X,Y) \leq 0$.*

Rules for the propagation of signs have been developed elsewhere [38, 39, 41] and, as seen from Proposition 2, are useful for determining the sign of covariances; however, as will be seen below, rules for deriving the sign of conditional covariances are more subtle. Theorem 4 concerns the conditional covariance of two parents of the node $D$. However, often what will be desired is the sign of the conditional covariance of two variables which are not parents of the conditioning node. For example, in the coaggregation problem discussed in the Introduction, we wanted to draw conclusions about $\mathrm{Cov}(P_2, B_1 | P_1 = 1)$, but neither $P_2$ nor $B_1$ are parents of $P_1$ in Figure 1. In the remainder of the paper we will thus extend Theorem 4 so as to allow for application to two variables, say $F$ and $G$, which are not parents of the conditioning node $D$. The variables $F$ and $G$ might be ancestors, descendants or have common causes with the parents, $E_1$ and $E_2$, of $D$. Consider, for example, the causal directed acyclic graphs in Figure 5.

If we were interested in the sign of $\mathrm{Cov}(F, G | D)$ in Figures 5(i)–(iii), then clearly Theorem 4 is insufficient. Theorems 5 and 6 below will allow us to extend the conclusions of Theorem 4 to examples such as those in Figure 5 and to certain other cases involving two variables that may not be parents of the conditioning variable. Lemmas 1–5 below will be needed in the proofs and application of Theorems 5 and 6. Lemmas 1 and 2 are consequences of Theorems 1 and 2 in the work of Esary, Proschan and Walkup [5]. Lemmas 3–5 are proved elsewhere in related work concerning the properties of monotonic effects [38].

LEMMA 1. *Let $f$ and $g$ be functions with $n$ real-valued arguments, such that both $f$ and $g$ are nondecreasing in each of their arguments. If $X = (X_1, \ldots, X_n)$ is a multivariate random variable with $n$ components, such that each component is independent of the other components, then $\mathrm{Cov}(f(X), g(X)) \geq 0$.*

LEMMA 2. *If $F$ and $G$ are binary and $u_1$ and $u_2$ are nondecreasing functions, then $\mathrm{sign}(\mathrm{Cov}(u_1(F), u_2(G))) = \mathrm{sign}(\mathrm{Cov}(F, G))$.*

LEMMA 3. *Let $X$ denote some set of nondescendants of $A$ that blocks all backdoor paths from $A$ to $Y$. If all directed paths between $A$ and $Y$ are positive, then $P(Y > y|a, x)$ and $\mathbb{E}[y|a, x]$ are nondecreasing in $a$.*

LEMMA 4. *Suppose that $E$ is binary. Let $Q$ be some set of variables which are not descendants of $F$ nor of $E$, and let $C$ be the common causes of $E$ and $F$ not in $Q$. If all directed paths from $E$ to $F$ (or from $F$ to $E$) are of positive sign and all directed paths from $C$ to $E$ not through $\{Q, F\}$ are of the same sign as all directed paths from $C$ to $F$ not through $\{Q, E\}$, then $E[F|E, Q]$ is nondecreasing in $E$.*

LEMMA 5. *Suppose that $E$ is not a descendant of $F$. Let $Q$ be some set of nondescendants of $E$ that block all backdoor paths from $E$ to $F$ and let $D$ be a node on a directed path from $E$ to $F$ such that all backdoor paths from $D$ to $F$ are blocked by $\{E, Q\}$. If all directed paths from $E$ to $F$, except possibly those through $D$, are of positive sign, then $\mathbb{E}[F|D, Q, E]$ is nondecreasing in $E$.*

Obvious analogues concerning negative signs hold for all of the lemmas above. Theorem 5 below will allow us to determine the sign of the conditional covariance of $F$ and $G$ on graphs like those in Figure 5, provided there are appropriate signs on the edges. The conclusion of Theorem 5 concerns the equality of the sign of two conditional covariances, $\mathrm{Cov}(F, G|D)$ and $\mathrm{Cov}(E_1, E_2|D)$. The theorem itself does not require knowledge of a sufficient causation representation and thus applies to general causal directed acyclic graphs. However, to draw conclusions about the sign of $\mathrm{Cov}(E_1, E_2|D)$, one must still appeal to Theorem 4 which does require some knowledge of a sufficient causation representation.

THEOREM 5. *Suppose that $E_1$, $E_2$ and $D$ are binary variables, that $E_1$ and $E_2$ are parents of $D$, that $F$ and $G$ are d-separated given $\{E_1, E_2, D\}$, that $F$ and $\{E_2, D\}$ are d-separated given $E_1$ and that $G$ and $\{E_1, D\}$ are d-separated given $E_2$. If $\mathrm{Cov}(F, E_1) \geq 0$ and $\mathrm{Cov}(G, E_2) \geq 0$ then $\mathrm{sign}(\mathrm{Cov}(F, G|D)) = \mathrm{sign}(\mathrm{Cov}(E_1, E_2|D))$.*

PROOF. Conditioning on $E_1$ and $E_2$, we have

$$\mathrm{Cov}(F, G|D) = \mathbb{E}[\mathrm{Cov}(F, G|D, E_1, E_2)|D] \\ + \mathrm{Cov}(\mathbb{E}[F|D, E_1, E_2], \mathbb{E}[G|D, E_1, E_2]|D).$$

The first expression is 0 since $F$ and $G$ are d-separated given $\{E_1, E_2, D\}$. Furthermore, since $F$ and $\{E_2, D\}$ are d-separated given $E_1$ and $G$ and $\{E_1, D\}$ are d-separated given $E_2$, the second expression can be reduced to $\mathrm{Cov}(\mathbb{E}[F|E_1], \mathbb{E}[G|E_2]|D)$. Thus,

$$\mathrm{Cov}(F, G|D) = \mathrm{Cov}(\mathbb{E}[F|E_1], \mathbb{E}[G|E_2]|D).$$

If $\mathrm{Cov}(F, E_1) \geq 0$ and $\mathrm{Cov}(G, E_2) \geq 0$ then, since $E_1$ and $E_2$ are binary, we have that $\mathbb{E}[F|E_1]$ is nonincreasing in $E_1$ and $\mathbb{E}[G|E_2]$ is nonincreasing in $E_2$, and so by Lemma 2, $\mathrm{sign}(\mathrm{Cov}(\mathbb{E}[F|E_1], \mathbb{E}[G|E_2]|D)) = \mathrm{sign}(\mathrm{Cov}(E_1, E_2|D))$. We thus have

$$\mathrm{sign}(\mathrm{Cov}(F, G|D)) = \mathrm{sign}(\mathrm{Cov}(E_1, E_2|D))$$

and this completes the proof. □

Note Theorem 5 requires that $\mathrm{Cov}(F, E_1) \geq 0$ and $\mathrm{Cov}(G, E_2) \geq 0$; Proposition 2 can be used to check whether these covariances are nonnegative; that is, the covariances will be nonnegative if $F$ and $E_1$ are positively monotonically associated and if $G$ and $E_2$ are positively monotonically associated.

EXAMPLE 5. Note that the graphs in Figures 5(i) and (ii) satisfy the d-separation restrictions of Theorem 5. In Figure 5(i), $G$ is an ancestor of $E_2$ whereas $F$ is related to $E_1$ as a descendant and by a common cause. In Figure 5(ii), $F$ is a descendant of $E_1$ and $G$ is related to $E_2$ both as an ancestor and by a common cause. The d-separation restrictions of Theorem 5 would still hold in Figures 5(i) and (ii) if $F$ and $E_1$ or $G$ and $E_2$ had multiple common causes or if there were several intermediate variables between $E_1$ and $F$ and between $G$ and $E_2$.

Note, however, that Theorem 5 requires that $F$ be $d$-separated from $\{E_2, D\}$ given $E_1$ and that $G$ be $d$-separated from $\{E_1, D\}$ given $E_2$. Thus, if $F$ or $G$ were a descendant of $D$, these assumptions would be violated. Consequently, Theorem 5 could not be applied to the diagram in Figure 5(iii). Nor could Theorem 5 be applied to the paper's introductory motivation to draw conclusions about the sign of $\mathrm{Cov}(P_2, B_1 | P_1 = 1)$ for the graph in Figure 1, since $B_1$ is a descendant of the conditioning variable $P_1$.

Theorem 6 below gives a result that allows for $F$ and $G$ to be descendants of $D$. Before stating this result we note, however, that Theorem 5 is restricted in yet another way. Theorem 5 required that $F$ and $G$ be $d$-separated given $\{E_1, E_2, D\}$. If $F$ and $G$ have common causes then the $d$-separation restrictions required by Theorem 5 will again, in general, not hold. Theorem 5 would thus not apply to the graphs given in Figure 6.

Theorem 6 gives a result similar to Theorem 5 which allows for $F$ or $G$ to be descendants of $D$ and allows also for $F$ and $G$ to have common causes. As with Theorem 5, the conclusion of Theorem 6 concerns the equality of the sign of two conditional covariances and the theorem itself does not require knowledge of a sufficient causation representation. But once again, to draw conclusions about the sign of $\mathrm{Cov}(F, G | D)$ using Theorem 6, one must know the sign of $\mathrm{Cov}(E_1, E_2 | D)$ and thus, appeal must again be made to Theorem 4 which does require some knowledge of a sufficient causation representation.

THEOREM 6. *Suppose that $E_1$, $E_2$ and $D$ are binary variables, that $E_1$ and $E_2$ are parents of $D$, that $F$ and $G$ are $d$-separated given $\{E_1, E_2, D, Q\}$, where $Q$ is some set of common causes of $F$ and $G$ (each component of which is univariate and independent of the other components in $Q$) that $F$ and $E_2$ are $d$-separated given $\{E_1, D, Q\}$, that $G$ and $E_1$ are $d$-separated given $\{E_2, Q, D\}$, that $Q$ and $\{E_1, E_2\}$ are $d$-separated given $D$ and that $Q$* and $D$ are $d$-separated. Suppose also that $\mathbb{E}[F|E_1, D, Q]$ is nondecreasing in $E_1$ and that $\mathbb{E}[G|E_2, D, Q]$ is nondecreasing in $E_2$. If $\mathrm{Cov}(E_1, E_2|D) \geq 0$, and for each element of $Q_i$ of $Q$, every directed path from $Q_i$ to $F$ is the same sign as every directed path from $Q_i$ to $G$, then $\mathrm{Cov}(F, G|D) \geq 0$. If $\mathrm{Cov}(E_1, E_2|D) \leq 0$, and for each element of $Q_i$ of $Q$, every directed path from $Q_i$ to $F$ is the opposite sign as every directed path from $Q_i$ to $G$, then $\mathrm{Cov}(F, G|D) \leq 0$.

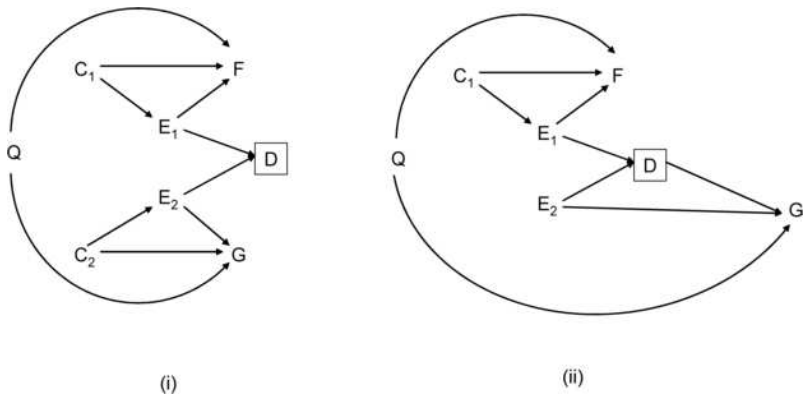

FIG. 6. *Examples in which $F$ and $G$ have a common cause.*

PROOF. We will prove the first of the results above; the proof of the second is similar. Conditioning on $\{E_1, E_2, Q\}$, we have

$$\mathrm{Cov}(F, G|D) = \mathbb{E}[\mathrm{Cov}(F, G|D, Q, E_1, E_2)|D] \\ + \mathrm{Cov}(\mathbb{E}[F|D, Q, E_1, E_2], \mathbb{E}[G|D, Q, E_1, E_2]|D).$$

The first expression is 0 since $F$ and $G$ are $d$-separated given $\{E_1, E_2, Q, D\}$. We can furthermore re-write the second expression as follows:

$$\mathrm{Cov}(F, G|D) \\ = \mathrm{Cov}(\mathbb{E}[F|D, Q, E_1, E_2], \mathbb{E}[G|D, Q, E_1, E_2]|D) \\ = \mathbb{E}[\mathrm{Cov}(\mathbb{E}[F|D, Q, E_1, E_2], \mathbb{E}[G|D, Q, E_1, E_2]|Q, D)|D] \\ + \mathrm{Cov}(\mathbb{E}[\mathbb{E}[F|D, Q, E_1, E_2]|Q, D], \mathbb{E}[\mathbb{E}[G|D, Q, E_1, E_2]|Q, D]|D).$$

We will show that each of these two expressions is positive. Since $F$ and $E_2$ are $d$-separated given $\{E_1, D, Q\}$, $\mathbb{E}[F|D, Q, E_1, E_2] = \mathbb{E}[F|E_1, D, Q]$; and since $G$ and $E_1$ are $d$-separated given $\{E_2, D, Q\}$, $\mathbb{E}[G|D, Q, E_1, E_2] = \mathbb{E}[G|E_2, D, Q]$. By assumption, we have that $\mathbb{E}[F|E_1, D, Q]$ is nondecreasing in $E_1$ and that $\mathbb{E}[G|E_2, D, Q]$ is nondecreasing in $E_2$. For fixed $q$,

$$\mathrm{Cov}(\mathbb{E}[F|D, Q=q, E_1, E_2], \mathbb{E}[G|D, Q=q, E_1, E_2]|Q=q, D) \\ = \mathrm{Cov}(\mathbb{E}[F|E_1, D, Q=q], \mathbb{E}[G|E_2, D, Q=q]|Q=q, D) \\ = \mathrm{Cov}(\mathbb{E}[F|E_1, D, Q=q], \mathbb{E}[G|E_2, D, Q=q]|D),$$

since $Q$ and $\{E_1, E_2\}$ are $d$-separated given $D$. And since $\mathbb{E}[F|E_1, D, Q=q]$ is nondecreasing in $E_1$ and $\mathbb{E}[G|E_2, D, Q=q]$ is nondecreasing in $E_2$, by Lemma 2, $\mathrm{Cov}(\mathbb{E}[F|E_1, D, Q=q], \mathbb{E}[G|E_2, D, Q=q]|D) = \mathrm{Cov}(E_1, E_2|D) \geq 0$. Thus, we have that $\mathrm{Cov}(\mathbb{E}[F|D, Q=q, E_1, E_2], \mathbb{E}[G|D, Q=q, E_1, E_2]|Q=q, D) \geq 0$ for all $q$ and taking expectations over $Q$ we have $\mathbb{E}[\mathrm{Cov}(\mathbb{E}[F|D, Q, E_1, E_2], \mathbb{E}[G|D, Q, E_1, E_2]|Q, D)|D] \geq 0$. We have shown that the first of the two expressions above is nonnegative. We now show that the second expression

$$\mathrm{Cov}(\mathbb{E}[\mathbb{E}[F|D, Q, E_1, E_2]|Q, D], \mathbb{E}[\mathbb{E}[G|D, Q, E_1, E_2]|Q, D]|D)$$

is also nonnegative. As before, $\mathbb{E}[F|D,Q,E_1,E_2] = \mathbb{E}[F|E_1,D,Q]$ and $\mathbb{E}[G|D, Q,E_1,E_2] = \mathbb{E}[G|E_2,D,Q]$. By hypothesis, for each element of $Q_i$ of $Q$ every directed path from $Q_i$ to $F$ is the same sign as every directed path from $Q_i$ to $G$; without loss of generality, we may assume that the sign of all of these directed paths are positive. By Lemma 3 with $X = \{E_1, D\}$ and $X = \{E_2, D\}$, respectively, $\mathbb{E}[F|E_1, D, Q = q]$ and $\mathbb{E}[G|E_2, D, Q = q]$ are both nondecreasing in each dimension of $q$. Note that we may apply Lemma 3 because if there were any backdoor paths from $Q$ to $F$ or to $G$, then $Q$ would have some parent which would also be a common cause of $F$ and $G$ and thus also a member of the set $Q$, but this would violate the assumption that the members of $Q$ were independent of one another. Furthermore,

$$\mathbb{E}[\mathbb{E}[F|D, Q = q, E_1, E_2]|Q = q, D] = \mathbb{E}[\mathbb{E}[F|E_1, D, Q = q]|Q = q, D]$$
$$= \mathbb{E}[\mathbb{E}[F|E_1, D, Q = q]|D]$$

and similarly, $\mathbb{E}[\mathbb{E}[G|D, Q = q, E_1, E_2]|Q = q, D] = \mathbb{E}[\mathbb{E}[G|E_2, Q = q]|D] = \mathbb{E}[\mathbb{E}[G|E_2, Q = q]|Q = q, D]$ since $Q$ and $\{E_1, E_2\}$ are $d$-separated given $D$. Thus,

$$\mathbb{E}[\mathbb{E}[F|D, Q = q, E_1, E_2]|Q = q, D] = \mathbb{E}[\mathbb{E}[F|E_1, D, Q = q]|D]$$

and

$$\mathbb{E}[\mathbb{E}[G|D, Q = q, E_1, E_2]|Q = q, D] = \mathbb{E}[\mathbb{E}[G|E_2, D, Q = q]|D]$$

are both nondecreasing in each dimension of $q$ from which it follows by Lemma 1 that $\text{Cov}(\mathbb{E}[\mathbb{E}[F|D, Q, E_1, E_2]|Q, D], \mathbb{E}[\mathbb{E}[G|D, Q, E_1, E_2]|Q, D]) \geq 0$. Since $Q$ and $D$ are $d$-separated we also have

$$\text{Cov}(\mathbb{E}[\mathbb{E}[F|D, Q, E_1, E_2]|Q, D], \mathbb{E}[\mathbb{E}[G|D, Q, E_1, E_2]|Q, D]|D)$$
$$= \text{Cov}(\mathbb{E}[\mathbb{E}[F|D, Q, E_1, E_2]|Q, D], \mathbb{E}[\mathbb{E}[G|D, Q, E_1, E_2]|Q, D]) \geq 0$$

and this completes the proof. □

Note the application of Theorem 6 requires that $\mathbb{E}[F|E_1, D, Q]$ is nondecreasing in $E_1$ and that $\mathbb{E}[G|E_2, D, Q]$ is nondecreasing in $E_2$. Either of the following will suffice for $\mathbb{E}[F|E_1, D, Q]$ to be nondecreasing in $E_1$ (similar remarks hold for $\mathbb{E}[G|E_2, D, Q]$): (i) $F$ and $D$ are $d$-separated given $\{Q, E_1\}$ and $F$ and $E_1$ are positively monotonically associated or (ii) if $F$ is a descendant of $E_1$ and $D$, $F$ and $E_1$ do not have common causes and all directed paths from $E_1$ to $F$ not through $D$ are of positive sign. Condition (i) suffices by Lemma 4; condition (ii) suffices by Lemma 5.

EXAMPLE 6. Although the graphs in Figure 5(iii) and in Figure 6 do not satisfy the $d$-separation restrictions of Theorem 5, it can be verified that the these graphs do satisfy the $d$-separation restrictions of Theorem 6.

At first glance, the $d$-separation restrictions of Theorems 5 and 6 appear to severely limit the settings to which conclusions about conditional covariances can be drawn. The $d$-separation requirements are, in fact, somewhat less restrictive than they may first seem. We argue that the $d$-separation restrictions of either Theorems 5 or 6 will apply to most graphs in which neither $F$ nor $G$ is a cause of the other (though the restrictions on the set of common causes $Q$, if any, of $F$ and $G$ in Theorem 6 are more substantial). Theorem 5 requires (i) that $F$ and $G$ are $d$-separated given $\{E_1, E_2, D\}$ and (ii) that $F$ and $\{E_2, D\}$ are $d$-separated given $E_1$ and that $G$ and $\{E_1, D\}$ are $d$-separated given $E_2$. In Theorems 5 and 6 (and Figures 5 and 6), $F$ was either an ancestor or descendant of or shared a common cause with $E_1$; and $G$ was either an ancestor or descendant of or shared a common cause with $E_2$. The $d$-separation restrictions essentially just require that $F$ and $G$ are sufficiently structurally separated so that (i) $F$ and $G$ are only associated because of $\{E_1, E_2, D\}$ and (ii) $F$ is associated with $\{E_2, D\}$ only through $E_1$; and $G$ is associated with $\{E_1, D\}$ only through $E_2$. If neither $F$ or $G$ is a descendant of $D$, then the conditions will, in general, only be violated if one of $F$ or $G$ is a cause of the other or if they share a common cause. Theorem 6, however, allowed for $F$ and $G$ to have common causes $Q$. The restrictions on $Q$ in Theorem 6 were somewhat substantial, but the restrictions on $F$ and $G$ are very similar to those of Theorem 5 except that they were made conditional on $Q$. Theorems 5 and 6 will thus apply to a wide range of graphs, as can also be seen by the variety of graphs in Figures 5 and 6, in which neither $F$ nor $G$ is a cause of the other.

As is clear from Proposition 2, rules concerning the propagation of signs were sufficient to determine the sign of the covariance between two variables. For conditional covariances, the principles guiding such a determination are more subtle. The principle behind the proofs of Theorems 5 and 6 was to partition the conditional covariance into two components

$$\mathrm{Cov}(F, G|D) = \mathbb{E}[\mathrm{Cov}(F, G|D, Q, E_1, E_2)|D]$$
$$+ \mathrm{Cov}(\mathbb{E}[F|D, Q, E_1, E_2], \mathbb{E}[G|D, Q, E_1, E_2]|D)$$

with $Q = \varnothing$ in the proof of Theorem 5. The $d$-separation restrictions allowed for the conclusion that $\mathrm{Cov}(F, G|D, Q, E_1, E_2) = 0$. Additional $d$-separation restrictions were needed so that the second expression $\mathrm{Cov}(\mathbb{E}[F|D, Q, E_1, E_2], \mathbb{E}[G|D, Q, E_1, E_2]|D)$ could be reduced to a form in which the sign of this conditional covariance could be determined from signed edges and an appeal to Theorem 4.

Having stated Theorem 6, we can now return to the motivating example presented in the paper's Introduction.

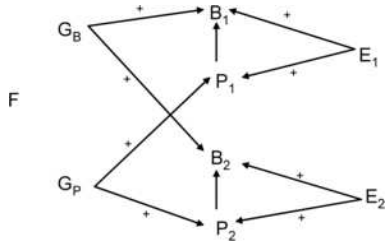

FIG. 7. *Causal directed acyclic graph with signed edges, under the null hypothesis of no familial coaggregation.*

EXAMPLE 7. In the motivating example described in Figure 1, with data available only on $P_1, P_2, B_1, B_2$, we wish to test the null hypothesis of no familial coaggregation (i.e., the null hypothesis that there are no directed edges emanating from $F$). Note that Hudson et al. [10] consider an alternative approach using a threshold model with additive multivariate normal latent factors. Here we use a sufficient causation approach. Given the substantive knowledge that for no subset of the population do the genetic causes $G_p$ and $G_B$ of $P$ and $B$ prevent disease and that for no subset of the population do the environmental causes $E_1$ and $E_2$ of $B$ and $P$ prevent either disease, we have that $E_1$ and $E_2$ have positive monotonic effects on $P_1$ and $B_1$ and on $P_2$ and $B_2$, respectively, and that $G_P$ has a positive monotonic effect on $P_1$ and on $P_2$ and that $G_B$ has a positive monotonic effect on $B_1$ and on $B_2$. The null hypothesis of no familial coaggregation can then be represented by the signed causal directed acyclic graph given in Figure 7.

If, in addition, using prior biological knowledge, it is assumed that there is no synergism between $E_1$ and $G_P$ in the sufficient cause sense, then we can apply part (vii) of Theorem 4 and, under the null hypothesis of no familial coaggregation, we have that $\text{Cov}(E_1, G_P | P_1 = 1) \leq 0$. By Theorem 6 with $Q = \varnothing$ we have that $\text{sign}(\text{Cov}(B_1, P_2 | P_1 = 1)) = \text{sign}(\text{Cov}(E_1, G_P | P_1 = 1))$. Under the null hypothesis of no familial coaggregation we thus have $\text{sign}(\text{Cov}(B_1, P_2 | P_1 = 1)) = \text{sign}(\text{Cov}(E_1, G_P | P_1 = 1)) \leq 0$. Thus, as claimed in the Introduction, a test of the null $\text{Cov}(B_1, P_2 | P_1 = 1) \leq 0$ is a test of no familial coaggregation under the assumption of no synergism between $E_1$ and $G_P$. Note that by the symmetry of this example, a test of the null $\text{Cov}(B_2, P_1 | P_2 = 1) \leq 0$ is a test of no familial coaggregation under the assumption of no synergism between $E_2$ and $G_P$. The development of a theory of minimal sufficient causation on directed acyclic graphs provided the concepts necessary to derive these results.

**6. Discussion.** In this paper we have incorporated notions of minimal sufficient causation into the directed acyclic graph causal framework. Doing so has provided a clear theoretical link between two major conceptualizations of causality. Causal directed acyclic graphs with minimal sufficient causation structures have furthermore allowed for the development of rules governing the sign of conditional covariances and of rules governing the presence of conditional independencies which hold only in a particular stratum of the conditioning variable.

The present work could be extended in a number of directions. Theory could be developed concerning cases in which a sufficient causation structure involves redundant sufficient causes or sufficient causes that are not minimally sufficient. Specifically, it might be possible to develop a system of axiomatic rules which govern conditional independencies within strata of variables on a causal directed acyclic graph with a sufficient causation structure, to furthermore demonstrate the soundness and completeness of this axiomatic system and to construct algorithms for applying the rules to identify all conditional independencies inherent in the graph's structure. Another direction of further research might involve the incorporation of the AND and OR nodes that arise from sufficient causation structures into other graphical models such as summary graphs [4], MC-graphs [12], chain graph models [2, 6, 14, 15, 16, 23, 34, 42] and ancestral graph models [24]. Finally, further work could be done extending the results of Theorem 4 to yet more general settings than those of Theorems 5 and 6.

DEPARTMENT OF HEALTH STUDIES
UNIVERSITY OF CHICAGO
5841 SOUTH MARYLAND AVENUE
MC 2007
CHICAGO, ILLINOIS 60637
USA
E-MAIL: vanderweele@uchicago.edu

HARVARD SCHOOL OF PUBLIC HEALTH
DEPARTMENTS OF BIOSTATISTICS
    AND EPIDEMIOLOGY
677 SOUTH HUNTINGTON AVENUE
BOSTON, MASSACHUSETTS 02115
USA
E-MAIL: robins@hsph.harvard.edu